\documentclass[12pt]{article}
\usepackage{amsmath,amssymb,amsfonts,amsthm}
\newcounter{item}[section]
\newcounter{kirshr}
\newcounter{kirsha}
\newcounter{kirshb}
\newenvironment{enumroman}{\setcounter{kirshr}{1}
\begin{list}{(\roman{kirshr})}{\usecounter{kirshr}} }{\end{list}}
\newenvironment{enumarab}{\setcounter{kirshb}{1}
\begin{list}{(\arabic{kirshb})}{\usecounter{kirshb}} }{\end{list}}
\newenvironment{athm}[1]{\vskip3mm\par\noindent
{\bf #1 }. \slshape }
{\upshape\par\vskip10pt minus3pt}
\newtheorem{theorem}{Theorem}[section]

\newtheorem{lemma}[theorem]{Lemma}

\newenvironment{demo}[1]{\noindent{\bf #1.}\upshape\mdseries}
{\nopagebreak{\hfill\rule{2mm}{2mm}\nopagebreak}\par\normalfont}
\theoremstyle{definition}

\newtheorem{definition}[theorem]{Definition}

\def\C{{\mathfrak{C}}}

\def\Nr{{\mathfrak{Nr}}}
\def\Fr{{\mathfrak{Fr}}}
\def\Sg{{\mathfrak{Sg}}}

\def\A{{\mathfrak{A}}}
\def\B{{\mathfrak{B}}}
\def\C{{\mathfrak{C}}}
\def\D{{\mathfrak{D}}}

\def\CA{{\bf CA}}

\def\SC{{\bf SC}}

\def\Dc{{\bf Dc}}

\def\(R)RA{{\bf (R)RA}}

\def\Dc{{\bf Dc}}

\def\Dc{{\bf Dc}}

 \def\CA{{\sf CA}}
\def\B{{\sf B}}

\def\Nr{{\mathfrak{Nr}}}

\def\Nr{{\mathfrak{Nr}}}

\def\A{{\mathfrak{A}}}
\def\B{{\mathfrak{B}}}
\def\C{{\mathfrak{C}}}
\def\D{{\mathfrak{D}}}

\def\CA{{\bf CA}}

\def\Nr{{\mathfrak{Nr}}}
\def\Fr{{\mathfrak{Fr}}}
\def\Sg{{\mathfrak{Sg}}}

\def\Ig{{\mathfrak{Ig}}}
\def\CA{{\bf CA}}

\def\SC{{\bf SC}}

\def\(R)RA{{\bf (R)RA}}

\def\Dc{{\bf Dc}}

\def\Dc{{\bf Dc}}

\def\MSA{{\bf MSA}}

\def\Nr{{\mathfrak{Nr}}}
\def\Fr{{\mathfrak{Fr}}}
\def\Sg{{\mathfrak{Sg}}}

\def\Ig{{\mathfrak{Ig}}}
\def\CA{{\bf CA}}

\def\SC{{\bf SC}}

\def\(R)RA{{\bf (R)RA}}

\def\Dc{{\bf Dc}}

\def\Dc{{\bf Dc}}

\def\Dc{{\bf Dc}}

 \def\CA{{\sf CA}}

\def\SC{{\bf SC}}

\def\A{{\mathfrak{A}}}
\def\B{{\mathfrak{B}}}
\def\C{{\mathfrak{C}}}
\def\D{{\mathfrak{D}}}

\def\Nr{{\mathfrak{Nr}}}
\def\F{{\mathfrak{F}}}
\def\CA{{\bf CA}}

\def\Sg{{\mathfrak Sg}}

\title{Interpolation for infinitary modal logic using polyadic algebras}
\author{Tarek Sayed Ahmed}

\begin{document}
\maketitle
 
\begin{abstract} We prove an interpolation theorem for certain infinitary modal logics in which atomic formulas are allowed to have infinite length. 

\end{abstract}

Modal cylindric  algebras were studied by Freeman \cite{Free} and modal polyadic algebras were studied by Georgescu \cite{G}, 
see \cite{HMT2} p.265-266 for an overview. The idea is to expand modal algebras by operations of cylindrifications and substitutions.
In the forementioned papers several quantified modal logics were proved (algebraically) to be complete with respect to usual 
Kripke semantics. In this paper we follow this trend, but our proven results are stronger in at least two respects. 
First we prove an interpolation theorem for several quantifed modal logics, of which the representability forementioned results readily follows,
and second our results apply to languages with an infinitary flavour; indeed atomic formulas can have infinite arity, though, like ordinary first 
order modal logic the scope of quantifiers 
is finite.

We prove an interpolation theorem for several quantified modal logics. We assume familiarity with the basic notions and concepts of (quantified )
modal logic \cite{b}, in particular with the modal operators of $\square$ (necessity) and $\Diamond$ (possibility), possible worlds, and Kripke 
semantics.  Propositional modal logics are defined in \cite{b} p. 36. We only deviate from this presentation by treating the implication 
symbol as a defined connective, not a basic one, and we treat the 
(possibility connective) $\Diamond$ as defined by $\neg \square\neg$, and not among
the primitive connectives.
The propositional modal logic $K$ (minimal modal logic) is defined on p.36 of \cite{b}. 
In the next page of op.cit the following propositional modal logics are defined:
$$D=K+\Diamond \top,$$
$$K4=K+(\square p\to \square\square p),$$ 
$$D4=D+(\square p \to \square\square p),$$ 
$$T=K+(\square p\to p),$$ 
$$S4=K4+(\square p\to p)$$ 
For each of the above calculi its relational semantics orginating from Hintikka and Kripke is considered.
$K$ is characterized by all Kripke models, $K4$ by models with transtive admissability relation and 
$S$ by reflexive and transitive models. When we are interested in the general study of a large family of logics, 
like all normal extensions of the minimal normal modal logic 
$K$, algebraic methods prove very fruitful.

We shall prove, using algebraic logic, that such logics with their quantified versions has the interpolation property. 
However, in our treatment of first order modal logic, we deviate from \cite{b} p.44 in allowing that atomic
formulas could be of infinite arity. Our results are therefore more general. 

\section{Algebraic Preliminaries}
All concepts and notions in this section is a simple modification of concepts already existing for cylindric algebras. All theorems
can be proved in a manner very similar to the cylindric case. Therefore we omit most proofs, 
instead we refer to their cylindric counterpart in \cite{HMT1}.
We provide sketches of proof whenever we find necessary.
\begin{definition} Let $\alpha$ be an ordinal. A modal substitution algebra of dimension $\alpha$, an $\MSA_{\alpha}$ for short, 
is an algebra of the form
$${\A}=(A, \lor, \land, 0,1, \square, {\mathsf c}_i, {\mathsf s}_i^j)_{i,j< \alpha}$$
where $(A, \lor, \land, -, 0, 1)$ is a Boolean algebra, $\square$, ${\mathsf c}_i, {\mathsf s}_i^j$ 
are unary operations on $\A$ ($i, j<\alpha$) satisfying the following 
equations for all $i,j,k,l < \alpha$:
\begin{enumerate}
\item $\square(-x\lor y)\leq (-\square x\lor -\square y).$
\item ${\mathsf c}_j0=0$
\item $x\leq {\mathsf c}_ix$
\item ${\mathsf c}_i (x\land {\mathsf c}_iy)={\mathsf c}_ix\land {\mathsf c}_iy$ 
\item $c_i(x\lor y)=c_ix\lor c_iy$
\item ${\mathsf c}_i{\mathsf c}_jx ={\mathsf c}_j{\mathsf c}_ix,$

\item ${\mathsf s}_i^ix=x,$

\item ${\mathsf s}_j^i$ are boolean endomorphisms,

\item  ${\mathsf s}_j^i{\mathsf c}_ix={\mathsf c}_ix,$

\item ${\mathsf c}_i{\mathsf s}_j^ix={\mathsf s}_j^ix$, whenever $i\neq j,$

\item  ${\mathsf s}_j^i{\mathsf c}_kx={\mathsf c}_k{\mathsf s}_j^ix$, whenever $k\notin \{i,j\},$

\item ${\mathsf c}_i{\mathsf s}_i^jx={\mathsf c}_j{\mathsf s}_j^ix,$



\item ${\mathsf s}_i^j{\mathsf s}_k^lx={\mathsf s}_k^l{\mathsf s}_i^jx$, whenever $|\{i,j,k,l\}|=4$


\item ${\mathsf s}_i^l{\mathsf s}_l^jx={\mathsf s}_i^l{\mathsf
s}_i^jx$.




\end{enumerate}
\end{definition}
We let ${\sf q_i}$ abbreviate $-{\sf c}_i-.$

The above is an algebraic reflection of quantified minimal logic $KQ$. Later, 
we shall impose more extra axioms concerning modalities, treating the quantified modal logics $DQ$, $TQ$, 
$KQ$, $D4Q$ and $S4Q$.

Examples of modal substitution algebras, arise from Kripke systems. 
For a relation $R$, we write  $xRy$ for $(x,y)\in R$.

\begin{definition} 
Let $W$ be a set and $R$ a binary relation of accessibility on $W$. Let $D$ be a non empty set (a domain). 
A Kripke system is a triple $\mathfrak{K}=( W,R, \{D_w\}_{w\in W})$ such that
for any $w\in W$, $D_w$ is a non-empty subset of $D$ satisfying 
$$wRw'\implies D_w\subseteq  D_{w'}.$$
\end{definition} 
Let $\mathfrak{O}$ be the Boolean algebra $\{0,1\}$. Now Kripke systems define concrete modal substitution algebras as follows.
Let $\alpha$ be any set and $\mathfrak{K}=(W, R, \{D_w\}_{w\in W})$ be a Kripke system.
Consider the set
$$\mathfrak{F}_{\mathfrak{K}}=\{(f_w:w\in W); f_w:{}^{\alpha}D_w\to \mathfrak{O}, wR w'\implies f_w\leq f_{w'}\}.$$
If $x,y\in {}^{\alpha}D_w$ and $j\in \alpha$ then we write $x\equiv _jy$ if $x(i)=y(i)$ for all $i\neq j$.
We write $(f_w)$ instead of $(f_w:w\in W)$.
In $\mathfrak{F}_{\mathfrak{K}}$ we consider the following operations:
$$(f_w)\lor (g_w)=(f_w\lor g_w)$$
$$(f_w)\land (g_w)=(f_w\land g_w.)$$
For any $(f_w)$ and $(g_w)\in \mathfrak{F}$, define
$$\neg (f_w)=(-f_w).$$
For any $(f_w)$ and $(g_w)\in \F$, define
$$\square(f_w)=(h_w),$$
where $(h_w)$ is given by $h_w(x)=1$ if and only if for any $w'\geq w$ $f_{w'}(x)=1$.
For any $i, j\in  {}^{\alpha}\alpha$, we define
$${\sf s}_i^j:\mathfrak{F}\to \mathfrak{F}$$by
$${\sf s}_i^j(f_w)=(g_w)$$
where
$$g_w(x)=f_w(x\circ [j,i])\text { for any }w\in W\text { and }x\in {}^{\alpha}D_w.$$
For any $j\in \alpha$ and $(f_w)\in \mathfrak{F}$
define
$${\sf c}_{j}(f_w)=(g_w)$$
where for $x\in {}^{\alpha}D_w$
$$g_w(x)=\bigvee\{f_w(y): y\in {}^{\alpha}D_w, y\equiv_jx\}.$$
It is tedious, but basically routine to verify  that
$\mathfrak{F}$ endowed with the above operations is an $\MSA_{\alpha}.$ 
We shall also deal with cases when the domain of $f_w$ is only a subset of $^{\alpha}D_w$, so that our semantics will be relativized.
\begin{definition}\begin{enumarab} 
\item Let $\A\in \MSA_{\alpha}$. Let $x\in A$. 
Then $\Delta x$, the dimension set of $x$, is the set
$\{i\in \alpha: {\sf c}_ix\neq x\}$. If we view $x$ as a formula, then 
$\Delta x$ represents the free variables occurring in $x.$

\item Let $\alpha\leq \beta$ be ordinals. Let ${\A}\in \MSA_{\beta}$. 
Then the  neat $\alpha$ reduct of $\A$ is the algebra with universe 
$Nr_{\alpha}A=\{x\in A: \Delta x\subseteq \alpha\}$
and whose operations are those of the similarity type of $\MSA_{\alpha}$ restricted to
$Nr_{\alpha}A$.  We denote such an algebra by $\Nr_{\alpha}\A$.
It is not hard to see that $\Nr_{\alpha}{\A}\in \MSA_{\alpha}$.

\item Let $L\in \MSA_{\beta}$. Then 
$Nr_{\alpha}L$ stands for all class of all neat $\alpha$ reducts in $L$.
\end{enumarab}
\end{definition}

Let $\alpha$ be an (infinite) ordinal. 
Then $\Dc_{\alpha}$ denotes the class of dimension 
complemented algebras in $\MSA_{\alpha}$.
$\A$ is such if $|\alpha\sim\Delta x|\geq \omega$ for all $x\in A$. 
This means that in languages we consider (atomic) formulas can be of infinite arity; the point is, 
like first order modal logic, there are infinitely many variable that do not occur in any given formula, but however, unlike first order modal logic,
the  variables occuring in a formula can well be infinite.

A $\Dc_{\alpha}$ free algebra cannot be dimension complemented \cite[2.5.15]{HMT2}.
This difficulty is easily 
conquered however by an appropriate relativisation of the concept of freedom, 
more precisely:

\begin{definition}{\cite[2.5.31]{HMT1}}
Let $\delta$ be a cardinal.
Let $\alpha$ be an ordinal. 
Let $_{\alpha} \Fr_{\delta}$ be the absolutely free algebra on $\delta$ 
generators and of type $\MSA_{\alpha}.$ 
For an algebra
$\cal A,$ we write $R\in Con\cal A$ if $R$ is a congruence relation on $\cal A.$
Let $\rho\in {}^{\delta}\wp(\alpha)$.  
Let $L$ be a class having the same similarity type as
$\MSA_{\alpha}.$ $SL$ denotes the class of all subalgebras of members of $L$.
Let
$$Cr_{\delta}^{(\rho)}L=\bigcap\{R: R\in Con_{\alpha}\Fr_{\delta}, 
{}_{\alpha}\Fr_{\delta}/R\in SL, 
c_k^{_{\alpha}\Fr_{\delta}}{\eta}/R=\eta/R \text { for each }$$
$$\eta<\delta \text 
{ and each }k\in \alpha\smallsetminus 
\rho{\eta}\}$$ 
and
$$\Fr_{\delta}^{\rho}L={}_{\alpha}\Fr_{\beta}/Cr_{\delta}^{(\rho)}L.$$
\end{definition}
The ordinal $\alpha$ does not appear in $Cr_{\delta}^{(\rho)}L$ and $\Fr_{\delta}^{(\rho)}L$
though it is involved in their definition. 
However, $\alpha$ will be clear from context so that no confusion is likely to ensue.
The algebra $\Fr_{\delta}^{(\rho)}L$ is referred to \cite{HMT1}
as a dimension restricted free algebra over $K$ with $\beta$
generators. Also $\Fr_{\delta}^{(\rho)}L$ is said to be 
dimension restricted by the function $\rho$, or simply, 
$\rho$-dimension-restricted.

\begin{definition} Let $\delta$ be a cardinal. 
Assume that $L\subseteq 
\MSA_{\alpha}$, $x=\langle x_{\eta}: \eta<\delta\rangle\in {}^{\delta}A$
and $\rho\in {}^{\delta}\wp(\alpha)$. Then we say that the sequence 
$x$ $L$-freely generates $\A$ under the 
dimension restricting function $\rho,$ if the following two conditions are satisfied:
\begin{enumroman}
\item ${\A}$ is generated by $Rgx$, 
and $\Delta x_{\eta}\subseteq \rho(\eta)$ for every $\eta<\delta.$ 

\item Whenever ${\B}\in L$, $y=\langle y_{\eta}: \eta<\delta\rangle\in   {}^{\delta}B$ 
and $\Delta y_{\eta}\subseteq \rho(\eta)$
for every $\eta <\delta,$ there is a homomorphism 
$h$ from $\A$ to $\B$
such that $h\circ x=y.$
\end{enumroman}
\end{definition}
For an algebra $\A$ and $X\subseteq A,$ we write, 
following \cite{HMT1}, $\Sg^{\A}X$, or even simply $\Sg X,$ 
for the subalgebra of $\A$
generated by $X$. We have the following:

\begin{theorem}\label{neat}
\begin{enumroman}

\item Let $L\subseteq \MSA_{\alpha}.$ Let $\delta$ be a cardinal. 
Let $\rho\in {}^{\delta}\wp(\alpha)$. 
Then the sequence $\langle \eta/Cr_{\delta}^{(\rho)}L: \eta<\delta\rangle$
$L$-freely generates $\Fr_{\delta}^{\rho}L$ 
under the dimension restricting function
$\rho$.

\item Let $\beta\geq \alpha$ be ordinals, $\delta$ be any cardinal $>0$
and $\rho\in {}^{\delta}\wp(\alpha)$. Let ${\B}\in \MSA_{\beta}$ and $x\in {}^{\delta}B$. 
If the sequence $x$ $\MSA_{\beta}$-freely generates
$\B$ under the dimension restricting function $\rho$, 
then $x$, $\Nr_{\alpha}\MSA_{\beta}$-freely generates
$\Sg^{\Nr_{\alpha}{\B}}Rgx$ under the 
dimension restricting function $\rho.$ 
\end{enumroman}
\end{theorem}
\begin{demo}{Proof} \cite{HMT1} theorem 2.6.45.
\end{demo}
\begin{definition} {\ }\label{d5-1}
\begin{enumroman}
\item
For a transformation $\tau\in {}^{\alpha}\alpha$, the support of
$\tau$, or $sup\tau$ for short, is the set $\{i\in \alpha:
\tau(i)\neq i\}$. $\tau$ is finite if $sup\tau$ is finite. For a
finite transformation $\tau$ we write $[u_0|v_0, u_1|v_1,\ldots,
u_{k-1}|v_{k-1}]$ if $sup\tau=\{u_0,\ldots ,u_{k-1}\}$, $u_0<u_1
\ldots <u_{k-1}$ and $\tau(u_i)=v_i$ for $i<k$.

\item Let $\A\in \Dc_{\alpha}$, that is $\alpha\sim \Delta x$ is
infinite for every $x\in A$. If $\tau=[u_0|v_0, u_1|v_1,\ldots,
u_{k-1}|v_{k-1}]$ is a finite transformation, if $x\in A$ and if
$\pi_0,\ldots ,\pi_{k-1}$ are in this order the first $k$ ordinals
in $\alpha\smallsetminus (\Delta x\cup Rg(u)\cup Rg(v))$, then
$${\mathsf s}_{\tau}x={\mathsf s}_{v_0}^{\pi_0}\ldots
{\mathsf s}_{v_{k-1}}^{\pi_{k-1}}{\mathsf s}_{\pi_0}^{u_0}\ldots
{\mathsf s}_{\pi_{k-1}}^{u_{k-1}}x.$$
\end{enumroman}
\end{definition}
We represent the function
$ {\mathsf s}^{\mu_{\pi_0}}_{\nu_{\pi_0}} \circ {\mathsf
s}^{\mu_{\pi_1}}_{\nu_{\pi_1}} \circ ......\circ {\mathsf
s}^{\mu_{\pi_{k-1}}}_{\nu_{\pi_{k-1}}} $ simply as ${\mathsf
s}^\mu_\nu$, whenever $\mu, \nu$ are functions from the same finite
subset of $ \omega $ into $\alpha$ and $\pi$ is the unique strictly
increasing sequence such that $ Rg (\pi) = Do (\mu) = Do (\nu) $.

\begin{lemma} Let $\A\in \MSA_{\alpha}$. The following statements
hold for all $x \in A$ and all $i,j,k,l<\alpha$:
\begin{enumarab}
\item ${\mathsf s}_j^i{\mathsf s}_i^k{\mathsf c}_ix={\mathsf s}_j^k{\mathsf c}_ix$.

\item $ {\mathsf s}^i_j {\mathsf s}^l_i {\mathsf c}_i {\mathsf c}_k x
= {\mathsf s}^k_j {\mathsf s}^l_k {\mathsf c}_i {\mathsf c}_k x.$

\item ${\mathsf s}^i_j {\mathsf s}^i_l x = {\mathsf s}^i_l x$,  if $i \neq l.$

\item ${\mathsf s}^i_j {\mathsf s}^j_i x = {\mathsf s}^i_j x.$

\item For all $i, j < \alpha$, $ {\mathsf s}_{[i|j]} x = {\mathsf s}^i_j x.$

Assume that in addition $\alpha\smallsetminus \Delta x$ is infinite
for every $x\in \A$, then for any such $x$ we have:

\item If $k < \omega$, if $ \mu, \nu, \pi, \rho \in  {}^k \alpha$, if
$ \mu, \pi, \rho$ are one-one, and if $ (Rg (\pi) \cup Rg (\rho))
\cap (\Delta x \cup Rg (\mu) \cup Rg (\nu)) = 0 $, then $ {\mathsf
s}^\pi_\nu {\mathsf s}^\mu_\pi x = {\mathsf s}^\rho_\nu {\mathsf
s}^\mu_\rho x.$
\item let $\tau$ be a finite transformation of $\alpha$. If $k <
\omega$, if $\mu, \pi \in {}^k \alpha $, if $\mu$ and $\pi$ are
one-one, if $\{ \lambda : \tau(\lambda) \neq \lambda \} \subseteq Rg
(\mu)$, $Rg(\pi)\cap (\Delta x\cup Rg(\mu)\cup \tau(Rg(\mu)))=0$,
then $ {\mathsf s}_\tau x = {\mathsf s}^\pi_{\tau \circ \mu}
{\mathsf s}^\mu_\pi x. $

\end{enumarab}
\end{lemma}

\begin{demo}{Proof} With respect to a different axiomatization, this is proved in \cite{HMT1} p.237. 
We prove $(1),(2)$ and $(6)$ to make sure that our axiomatization works. (The rest is similar, by using \cite{HMT1}.)
 \begin{equation*}
\begin{split} (1) \quad \qquad \qquad
{\mathsf s}^i_j {\mathsf s}^k_i {\mathsf c}_i x &= {\mathsf s}^i_j
{\mathsf s}^k_j
{\mathsf c}_i x \\
 &= {\mathsf s}^i_j {\mathsf c}_i {\mathsf s}^k_j  x \\
 &=  {\mathsf c}_i {\mathsf s}^k_j  x \\
&= {\mathsf s}^k_j  {\mathsf c}_i  x .\\
\end{split}
\end{equation*}

 \begin{equation*}
\begin{split} (2)   \quad \qquad \qquad
{\mathsf s}^i_j {\mathsf s}^l_i {\mathsf c}_i {\mathsf c}_k x &=
{\mathsf s}^l_i {\mathsf c}_i {\mathsf c}_k
 x \\
 &= {\mathsf s}^l_i {\mathsf c}_k {\mathsf c}_i x\\
 &=   {\mathsf s}^k_j {\mathsf s}^l_k {\mathsf c}_k {\mathsf c}_i x \\
&= {\mathsf s}^k_j {\mathsf s}^l_k {\mathsf c}_i {\mathsf c}_k x\\
\end{split}
\end{equation*}

(6) can be established if $Rg (\pi) \cap  Rg (\rho) = 0$ as follows
\begin{equation*}
\begin{split}
{\mathsf s}^\pi_\nu {\mathsf s}^\mu_\pi x &= {\mathsf
s}_{\nu_0}^{\pi_0}\ldots {\mathsf s}_{\nu_{k-1}}^{\pi_{k-1}}{\mathsf
s}_{\pi_0}^{\mu_0}\ldots {\mathsf s}_{\pi_{k-1}}^{\mu_{k-1}}x\\
 &= {\mathsf s}_{\nu_0}^{\pi_0}{\mathsf s}_{\pi_0}^{\mu_0} \ldots
 {\mathsf s}_{\nu_{k-1}}^{\pi_{k-1}}{\mathsf s}_{\pi_{k-1}}^{\mu_{k-1}}x\\
&= {\mathsf s}_{\nu_0}^{\rho_0}{\mathsf s}_{\rho_0}^{\mu_0}{\mathsf
c}_{\pi_0}
  {\mathsf c}_{\rho_0} \ldots
 {\mathsf s}_{\nu_{k-1}}^{\rho_{k-1}}{\mathsf s}_{\rho_{k-1}}^{\mu_{k-1}}
 {\mathsf c}_{\pi_{k-1}}
  {\mathsf c}_{\rho_{k-1}} x \\
&= {\mathsf s}_{\nu_0}^{\rho_0}{\mathsf s}_{\rho_0}^{\mu_0} \ldots
 {\mathsf s}_{\nu_{k-1}}^{\rho_{k-1}}{\mathsf s}_{\rho_{k-1}}^{\mu_{k-1}}
  x \\
&= {\mathsf s}_{\nu_0}^{\rho_0}\ldots {\mathsf
s}_{\nu_{k-1}}^{\rho_{k-1}}{\mathsf s}_{\rho_0}^{\mu_0}\ldots
{\mathsf s}_{\rho_{k-1}}^{\mu_{k-1}}x \\
&= {\mathsf s}^\rho_\nu {\mathsf s}^\mu_\rho x.
\end{split}
\end{equation*}

To prove (6) in the general case, note that there is a $ \xi \in
{}^k \alpha$ with $\xi$ one-one, $Rg (\xi) \cap (\Delta x \cup Rg
(\mu) \cup Rg (\nu)) = 0,$ and $Rg (\xi) \cap (Rg (\pi) \cup Rg (
\rho)) = 0, $ and then use the special case just established.
\end{demo}

\begin{theorem}\label{T5-4} Let $\A\in \SC_{\alpha}$ such that $\alpha\smallsetminus \Delta x$
is infinite for every $x$. For finite $\Gamma\subseteq \alpha$,
${\mathsf c}_{(\Gamma)}x$ abbreviates ${\mathsf c}_{i_0}\ldots
{\mathsf c}_{i_n}x$, where $\Gamma=\{i_0,\ldots, i_{n}\}$. Let
$\sigma$ and $\tau$ be finite transformations on $\alpha$, $\Gamma,
\Delta$ finite subsets of $\alpha$ and $x\in A$.

Then
\begin{enumroman}
\item ${\mathsf s}_{\tau}$ is a boolean endomorphism of $\A$.
\item ${\mathsf s}_{\sigma\circ \tau}={\mathsf s}_{\sigma}\circ {\mathsf s}_{\tau}$.
\item If $\sigma \upharpoonright (\alpha\smallsetminus \Gamma) = \tau \upharpoonright (\alpha \smallsetminus
 \Gamma )$, then ${\mathsf s}_{\sigma}{\mathsf c}_{(\Gamma)} x = {\mathsf s}_{\tau}{\mathsf c}_{(\Gamma)}x.$
\item if $\sigma\upharpoonright \Delta x=\tau\upharpoonright \Delta x$,
then ${\mathsf s}_{\sigma}x={\mathsf s}_{\tau}x$.
\item If $\tau^{-1}\Gamma=\Delta$ and $\tau\upharpoonright \Delta $ is
one to one, then ${\mathsf c}_{(\Gamma)}{\mathsf s}_{\tau}x={\mathsf
s}_{\tau}{\mathsf c}_{(\Delta )}x$.
\end{enumroman}
\end{theorem}

\begin{demo}{Proof} Proceed like \cite{HMT1} p. 238, using the previous lemma.
\end{demo}
\begin{theorem} Let $\alpha<\beta$ be infinite ordinals. 
Then if $\A\in \Dc_{\alpha}$, then there exists $\B\in \Dc_{\beta}$ such that $\A\subseteq \Nr_{\alpha}\B$, and for all $X\subseteq A$,
one has $\Sg^{\A}X=\Nr_{\alpha}\Sg^{\B}X.$ Furthermore $\B$ is unique up to an isomorphism that fixes $A$ pointwise.
\end{theorem} 
\begin{demo}{Proof} The proof is completely analogous to \cite{HMT1}, theorems 2.6.49, 2.6.67 and 2.6.72.
\end{demo}
We call $\B$ in the previous theorem the $\beta$ dilation of $\A.$
\begin{athm}{Lemma 2.5} Assume $\omega\leq \alpha\leq \beta$, and $\delta$ is a cardinal.
Let $\rho\in {}^{\delta }\wp(\alpha)$ be such that 
$|\alpha \smallsetminus \rho(i)|\geq \omega$ for all $i\in \delta$. 
Then 
$\Fr_{\delta}^{\rho}\MSA_{\alpha}\cong \Nr_{\alpha}\Fr_{\delta}^{\rho}\MSA_{\beta}.$
\end{athm}
\begin{demo}{Proof}
Let ${\A}\in \MSA_{\alpha}$ and $a=\langle a_{\eta} :\eta<\delta \rangle\in {}^{\delta}A$
be such that $\rho(\eta)$ includes the dimension set of $a_{\eta}$ for every $\eta<\delta.$
Since $\Fr_{\delta}\MSA_{\alpha}\in \Dc_{\alpha}$, we have by \cite[2.1.7]{HMT1}
and \cite[2.5.33(ii)]{HMT1} that 
$|\alpha\smallsetminus \cup_{\eta\in \Gamma}\rho(\eta)|\geq \omega$ for each finite
$\Gamma\subseteq \delta$.  We can assume without loss of generality that 
$Rga$ generates $\A$. Then by \cite[2.1.7]{HMT1} we see that 
${\A}\in \Dc_{\alpha}$. Thus by \cite [2.6.34, 2.6.35, 2.6.49]
{HMT1}, we have ${\A}\in Nr_{\alpha}\MSA_{\beta}.$
Let ${\D}=\Fr_{\delta}^{\rho}\MSA_{\beta}$. Then  by theorem \ref{neat}, we have 
$x=\langle \eta/Cr_{\delta}^{(\rho)}\MSA_{\beta}:  \eta<\delta\rangle$ 
$Nr_{\alpha}\CA_{\beta}$-freely generates $\Sg^{\Nr_{\alpha}\D}Rgx.$
Therefore there exists a homomorphism 
$h:\Sg^{\Nr_{\alpha}\D}Rgx \to \A$
such that $h(\eta/Cr_{\delta}^{\rho}\MSA_{\beta})=a_{\eta}$ 
for all $\eta<\delta$. 
But by \cite[ 2.6.67 (ii)]{HMT1} we have
$$\Nr_{\alpha}\Fr_{\delta}^{\rho}\MSA_{\beta}=
\Nr_{\alpha}\Sg^{\D}Rgx=\Sg^{\Nr_{\alpha}\D}Rgx$$
Therefore $\langle \eta/Cr_{\delta}^{\rho}\MSA_{\beta}: \eta<\delta\rangle$  
$\MSA_{\alpha}$-freely generates
$\Nr_{\alpha}\Fr_{\delta}^{\rho}\MSA_{\beta}.$
\end{demo}

Now we define our algebras. Their similarity type depends on a fixed in advance semigroup. 
We write $X\subseteq_{\omega} Y$ to denote that $X$ is a finite subset 
of $Y$.
\begin{definition} Let $\alpha$ be an infinite set. Let $G\subseteq {}^{\alpha}\alpha$ be a semigroup under the operation of composition of maps. 
An $\alpha$ dimensional polyadic Heyting $G$ algebra, a $GPHA_{\alpha}$ for short, is an algebra of the following
type
$$(A,\lor,\land,\rightarrow, 0, {\sf s}_{\tau}, {\sf c}_{(J)}, {\sf q}_{(J)})_{\tau\in G\subseteq {}^{\alpha}\alpha, J\subseteq_{\omega} \alpha}$$
where
$(A,\lor,\land, \rightarrow, 0)$ is a Heyting algebra, ${\sf s}_{\tau}:\A\to \A$ is a endomorphism of Heyting algebras,
${\sf c}_{(J)}$ is an existential quantifier, ${\sf q}_{(J)}$ is a universal quantifier, such that the following hold for all 
$p\in A$, $\sigma, \tau\in [G]$ and $J,J'\subseteq_{\omega} \alpha:$
\begin{enumarab}
\item ${\sf s}_{Id}p=p.$
\item ${\sf s}_{\sigma\circ \tau}p={\sf s}_{\sigma}{\sf s}_{\tau}p$ (so that $S:\tau\mapsto {\sf s}_{\tau}$ defines a homomorphism from $G$ to $End(\A)$; 
that is $(A, \lor, \land, \to, 0, G, S)$ is a transformation system).
\item ${\sf c}_{(J\cup J')}p={\sf c}_{(J)}{\sf c}_{(J')}p , \ \  {\sf q}_{(J\cup J')}p={\sf q}_{(J)}{\sf c}_{(J')}p.$
\item ${\sf c}_{(J)}{\sf q}_{(J)}p={\sf q}_{(J)}p , \ \  {\sf q}_{(J)}{\sf c}_{(J)}p={\sf c}_{(J)}p.$
\item If $\sigma\upharpoonright \alpha\sim J=\tau\upharpoonright \alpha\sim J$, then
${\sf s}_{\sigma}{\sf c}_{(J)}p={\sf s}_{\tau}{\sf c}_{(J)}p$ and ${\sf s}_{\sigma}{\sf q}_{(J)}p={\sf s}_{\tau}{\sf q}_{(J)}p.$
\item If $\sigma\upharpoonright \sigma^{-1}(J)$ is injective, then
${\sf c}_{(J)}{\sf s}_{\sigma}p={\sf s}_{\sigma}{\sf c}_{\sigma^{-1}(J)}p$
and ${\sf q}_{(J)}{\sf s}_{\sigma}p={\sf s}_{\sigma}{\sf q}_{\sigma^{-1}(J)}p.$
\end{enumarab}
\end{definition}

\begin{athm}{Notation.} For  a set $X$, recall that $|X|$ stands for the cardinality of $X$.
$Id_X$, or simply $Id$ when $X$ is clear from context,  denotes the identity function
on $X$.
For functions $f$ and $g$ and a set $H$, $f[H|g]$ is the function
that agrees with $g$ on $H$, and is otherwise equal to $f$. $Rng(f)$ denotes 
the range of $f$. For a transformation $\tau$ on $\omega$, 
the support of $\tau$, or $sup(\tau)$ for short, is the set:
$$sup(\tau)=\{i\in \omega: \tau(i)\neq i\}.$$
Let $i,j\in \omega$, then $\tau[i|j]$ is the transformation on $\omega$ defined as follows:
$$\tau[i|j](x)=\tau(x)\text { if } x\neq i \text { and }\tau[i|j](i)=j.$$ 
For a function $f$, $f^n$ denotes the composition $f\circ f\ldots \circ f$
$n$ times.
\end{athm}

\begin{definition}  (a) Let $\alpha$ be a countable set. Let $T\subseteq \langle {}^{\alpha}\alpha, \circ \rangle$ be a semigroup.
We say that $T$ is {\it rich } if $T$ satisfies the following conditions:
\begin{enumerate}
\item $(\forall i,j\in \alpha)(\forall \tau\in T) \tau[i|j]\in T.$
\item There exists $\sigma,\pi\in T$ such that
$(\pi\circ \sigma=Id,\  Rg(\sigma)\neq \alpha).$
\item $ (\forall \tau\in T)(\sigma\circ \tau\circ \pi)[(\alpha\sim Rg(\sigma))|Id]\in T.$
\item Let $T\subseteq \langle {}^{\alpha}\alpha, \circ\rangle$ be a rich semigroup. 
Let $\sigma$ and $\pi$ be as in (a) above. 
If $\sigma$ and $\pi$ satisfy $(i)$, $(ii)$ below:
$$ (i)\ \ \ \  (\forall n\in \alpha) |supp(\sigma^n\circ \pi^n)|<\alpha.$$  
$$(ii) \ \ \ \  (\forall n\in \alpha)[supp(\sigma^n\circ \pi^n)\subseteq 
\alpha\smallsetminus Rng(\sigma^n)];$$
then we say that $T$ is  {\it a strongly rich} semigroup. 
\end{enumerate}
\end{definition}

An example of a rich semigroup in $(^{\omega}\omega, \circ)$ and its  semigroup generated by
$\{[i|j], [i,j], suc, pred\}$. Here $suc$ abbreviates the successor function on $\omega$ and $pred$ is the function defined by 
$pred(0)=0$ and for other $n\in \omega$, $pred(n)=n-1$. In fact, both semigroups are strongly rich, 
in the second case $suc$ plays the role of $\sigma$ while $pred$ plays the role of 
$\pi$.

\section{Proof}
\begin{definition} Let $\A\in \Dc_{\alpha}$. Let $\Gamma,\Delta\subseteq \A$. 
We write $\Gamma \to _{\A}\Delta $, if there exist $a_1,\ldots a_n$ in $\Gamma $ and $b_1,\ldots b_m$ in $\Delta $ such that  
$a_1\land \ldots a_n\to b_1\lor\ldots b_m.$ We write $\Gamma\to_{\A}a$ if $\Delta=\{a\}$. By a
theory of  $\A\in \Dc_{\alpha},$ we mean a subset $T$ of $\A$
satisfying the condition: $(T\to _{\A}a)\Rightarrow a\in T$. If $F\subseteq A$ satisfies the dual condition:
$(a\rightarrow _{\A}F$) $\Rightarrow a\in F$, then $F$
is called a cotheory of $\A$. Denote by $Th(\A)$
the set of all theories and by $CTh(\A)$ the set of all cotheories of $\A$, respectively.
\end{definition}
\begin{definition} Let $\A\in \Dc_{\alpha}$ and $X_1, X_2\subseteq A$. A pair 
$(\Gamma ,\Delta )\in \Sg^{\A}X_1\times  \Sg^{\A}X_2$ is called $(\A, X_1, X_2)$ separable, or just $\A$ separable, 
whenever there exists $a\in \Sg^{\A}(X_1\cap X_2)$ such that 
$\Gamma \rightarrow _{\A}a$ and 
$a\rightarrow _{\A}\Delta $. Call a pair $(\A, X_1, X_2)$ inseparable, or simply inseparable, if it is not $\A$ separable. 
Note, that the pair $(T,F)$ is $\A$ inseparable if and only if 
$T\cap F=\emptyset $.

\end{definition}

\begin{definition}
Let $\A\in \Dc_{\alpha}$ and $X_1, X_2\subseteq A$.
A pair $(T,F)\in \A\times \A$ is called $(\A, X_1, X_2)$, or even simply $\A$ saturated , cf.
\cite{b} p. 134, if the following conditions are satisfied for all $a,b\in A$:

\begin{description}
\item  (1) $T\in Th(\Sg^{\A}X_1) $ and $F\in CTh(\Sg^{\A}X_2)$,

\item  (2) $T\cap F=\emptyset $, that is, $(T,F)$ is $\A$-inseparable,

\item  (3) $a\lor b\in T\Rightarrow (a\in T$ or $b\in T)$,

\item  (4) ${\sf c}_ka\in T\Rightarrow (\exists j\notin \Delta a)({\sf s}_j^ka\in T)$,

\item  (5) $a\land b\in F\Rightarrow (a\in F$ or $b\in F)$,

\item  (6) ${\sf  q}_ka\in F\Rightarrow (\exists j\notin \Delta a)({\sf s}_j^ka\in F)$.
\end{description}
\end{definition}

\begin{lemma} \label{neat} Let $\A\in \Dc_{\alpha}$. Let $X_1, X_2\subseteq A$. Let $(\Gamma, \Delta)$ be $(X_1, X_2, \A)$ inseparable. 
Then there exists an $\omega$ dilation $\B$ of $\A$ and $T\supseteq \Gamma$,
$F\supseteq \Delta$ such that $(T,F)$ is $(\B, X_1, X_2)$ saturated. 
\end{lemma}
\begin{demo}{Proof} \cite{b} lemma 5.6. Let $\A=\Nr_{\alpha}\B$,  with $\B$ the $\alpha+\omega$ dilation of $\A$. 
Let $\beta=\alpha+\omega$ be the dimension of 
$\B$.
Then $(\Gamma,\Delta)$ is $(\B,X_1,X_2)$ inseparable. If not, then there exists $c\in \Sg^{\B}(X_1\cap X_2)$ 
such that $\Gamma\to_{\B} c$ and $c\to_{\B} \Delta$. Now $|\Delta c\sim \alpha|<\omega$, since $A$ generates $\B$, hence 
there exists a finite subset $JJ$ of $\beta$ such that 
$$d={\sf c}_{(J)}c\in \Nr_{\alpha}\Sg^{\B}(X_1\cap X_2)=\Sg^{\Nr_{\alpha}\B}(X_1\cap X_2)=\Sg^{\A}(X_1\cap X_2),$$ and  
$\Gamma\to_{\A} d$ and $d\to_{\A} \Delta$ 
which is a contradiction.

Let $(a_1, \dots a_n\ldots)$ be an enumeration of $\Sg^{\B}X_1$ and $(b_1\ldots b_n,\ldots)$ be an enumeration of  $\Sg^{\B}X_2.$ 
This is possible, since $\A$ is countable, $A$ generates $\B$, and so $\B$ is countable.

We define $(T_n, F_n)$, $n\in \omega$ inductively. Let $T_0=\Gamma$ and $F_0=\Delta$. We assume inductively that 
$\beta\sim \bigcup_{x\in T_n}\cup\bigcup_{x\in F_n}\Delta x$
is infinite. This condition is clearly satisfied for the base of the induction, since $T_0\subseteq \A$ and $F_0\subseteq \A$, 
and $\B$ is the $\alpha+\omega$ 
dilation of $\A$.

At step $2n+1$. 

Let $F_{2n+1}=F_{2n}$ and to define $T_{2n+1}$ we distinguish between two cases:
 
{\bf Case 1}: The pair $(T_{2n}\cup \{a_n\}, F_{2n})$ is $\B$ separable. Then set $T_{2n+1}=T_{2n}.$
 
{\bf Case 2}: The pair $ (T_{2n}\cup \{a_{n}\}, F_{2n})$ is $\B$ inseparable.
We distinguish between two subcases:

{\bf Subcase 1}:  $a_{n}\neq {\sf c}_kx$, for any $k\in \beta$  and any $x\in B.$

Set   $T_{2n+1}=T_{2n}\cup \{a_{n}\}.$

{\bf Subcase 2}: $a_{n}={\sf c}_k x$, for some $k<\beta$ and some $x\in B.$ 

Let $T_{2n+1}=T_{2n}\cup \{a_{n}, {\sf s}_i^ka_n\}$ where  $i\notin \Delta a_n\cup\bigcup_{x\in T_n}\cup\bigcup_{x\in F_n}\Delta x.$ Such an $i$ exists 
by the induction hypothesis. 

Similarly at step $2n+2$, let $T_{2n+2}=T_{2n+1}$ and define $F_{2n+2}$ as follows: 
$F_{2n+2}=F_{2n+1}$, whenever the pair $(T_{2n+1},F_{2n+1}\cup \{b_{n}\})$ is $\B$-separable;
$F_{2n+2}=F_{2n+1}\cup \{b_{n},{\sf s}_i^jb_n\}$, whenever $b_{n}={\sf q}_i x$ for some $i<\beta$, $x\in B$ 
and the pair  $(T_{2n+1},F_{2n+1}\cup \{b_{n}\})$ is $\B$ inseparable, where 
$j\notin \Delta a_n\cup\bigcup_{x\in T_{2n+1}}\cup\bigcup_{x\in F_{2n+1}}\Delta x.$
One can prove by an easy induction on $n$ that each pair $(T_n, F_n)$ is $\B$ inseparable.
Now define 
\[
T=\{a\in \Sg^{\B}X_1:\exists n(T_{n}\rightarrow _{\B}a)\},F=\{b\in \Sg^{\B}X_2: \exists n(b\rightarrow _{\B}F_{n})\} 
\]
We check that $(T,F)$ is as required.  Clearly $T\in Th(\Sg^{\B}X_1)$ and $F\in CTh(\Sg^{\B}X_2)$.
One can check that $(T,F)$ is saturated by first observing that 
$a\in T$ iff $(T\cup \{a\}, F)$ is $\B$ inseparable and $b\in F$ iff $(T, F\cup \{b\})$ is $\B$ inseparable. Indeed, assume that $a=a_n$. 
The inseparability of $(T\cup \{a\}, F)$ implies the inseprability of $(T_{2n}\cup \{a\}), F_{2n})$ so $a\in T_{2n+1}\subseteq T$. 
On the other hand, from $a\in T$ 
the inseparability of $(T\cup \{a\}, F)$ follows since $(T,F)$ is inseparable.

We checkconditions  (3) and (4) of saturation. Assume that $a\lor b\in T$ but $a\notin T$ and  $b\notin T$. Then $(T\cup \{a\}, F)$ 
and $(T\cup \{b\}, F)$ would be separable. Hence there exists
$c\in \Sg^{\B}(X_1\cap X_2)$ such that $T\cup \{a\}\to_{\B} c$ and $c\to_{\B} F$ and there exists $c'\in \Sg^{\B}(X_1\cap X_2)$ 
such that $T\cup \{b\}\to_{\B} c'$ and $c'\to_{\B} F$, 
Then $T\cup \{a\lor b\}\to_{\B} c\lor c',$ $c\lor c'\to_{\B} F$ and $c\lor c'\in \Sg^{\B}(X_1\cap X_2)$. But this contradicts the inseparability of 
$(T,F)$.
Now we check (4). Assume that $a_n={\sf c}_ka\in T$ where $k<\beta$ and $a\in \B$. Then $(T_{2n}\cup \{a\}, F_{2n})$ is separable, 
so by construction ${\sf s}_k^ja\in T_{2n+1}\subseteq T.$
The rest of the conditions can be checked analogously.
\end{demo}
\begin{theorem} Let $\A=\Fr_{\beta}^{\rho}\Dc_{\alpha}$. Let $X_1, X_2\subseteq \A$. Assume that $(\Gamma, \Delta)$ is $(X_1, X_2, \A)$ 
inseparable. Then there is a Kripke system $\mathfrak{K}=(W, R, D_{w})_{w\in W}$ and a homomorphism $\psi:\A\to \F_{\mathfrak{K}},$ 
$w\in W,$ and $x\in {}^{\alpha}D_w$ such that  
$\psi(a)_w(x)=1$ for all $a\in \Gamma$ and $\psi(a)_w(x)=0$ for all $a\in \Delta.$
\end{theorem}

\begin{demo}{Proof} Let $\B_n=\Fr_{\beta}^{\rho}\Dc_{\alpha+n\omega}$, $n\leq \omega$. Then 
$$\A=\B_0\subseteq \B_1\subseteq\ldots \B_{\infty}=\Fr_{\beta}^{\rho}\Dc_{\alpha+\omega.\omega}$$
is a a sequence of minimal dilations. That is for $k<l\leq \omega$, we have $\B_k=\Nr_k\B_l$ and for all $X\subseteq B_k$, we have 
$\Sg^{\B_k}X=\Nr_k\Sg^{\B_l}X$.
Define inductively $W_n\subseteq \Sg^{\B_n}X_1\times \Sg^{\B_n}X_2$ such that $|W_n|\leq \omega$ and 
for $(T,F)\in W_n$, $(T, F)$ is $(X_1, X_2, \B_n)$ saturated   
and (*)
$$(\forall (T, F)\in W_n)( \forall a\in \Sg^{\B_n}X_1)(\Diamond a\in T\Rightarrow \exists ((T',F')\in W_{n+1})$$
$$[ \{a\}\cup \{a':\square a'\in \Gamma\}\subseteq T', \{b':\Diamond b'\in \Delta\}\subseteq F'].)$$
and
$$(\forall (T,F)\in W_n)( \forall b\in \Sg^{\B_n}X_2)(\square b\in F\Rightarrow \exists ((T',F')\in W_{n+1})$$
$$[ \{a':\square a\in \Gamma\}\subseteq T', \{\{b\}\cup \{b':\Diamond b'\in \Delta\}\subseteq F'].)$$
The base of the induction is provided by the previous lemma. 
Set $W_0=\{(T_0,F_0)\}$ where $(T_0, F_0)$ is a saturated extension of $(\Gamma ,\Delta)$ in $\B_1$. We have countably many tasks, 
so this can be accomplished once we make sure of:
\begin{athm}{Claim}\cite{b} lemma 5.5. Let $(\Gamma ,\Delta )$ be $(X_1, X_2, \B_n)$-inseparable. Then
\begin{description}
\item  (i) If $\Diamond a\in \Gamma $, then the pair $(\{a\}\cup \{a^{\prime
}|\square a^{\prime }\in \Gamma \},\{b^{\prime }|\Diamond b^{\prime }\in
\Delta \})$ is $(X_1, X_2, \B_n)$-inseparable.
\item  (ii) If $\square b\in \Delta $, then the pair $(\{a^{\prime }|\square
a^{\prime }\in \Gamma \},\{b\}\cup \{b^{\prime }|\Diamond b^{\prime }\in
\Delta \})$ is $(X_1, X_2,\B_n)$ inseparable.
\end{description}
\end{athm}
\noindent \textbf{Proof}
\begin{description} 
\item  (i) Assume that $((a\land a_1\land \ldots a_n\to_{\B_n} c)$ and $(c\to_{\B_n} (b_{1}\lor\ldots b_{m}))$ for some $c\in \Sg^{\B}(X_1\cap X_2)$, 
where $\square a_{1},\ldots ,\square a_{n}\in \Gamma $ 
and $\Diamond b_{1},\ldots ,\Diamond
b_{m}\in \Delta $. Then $\Diamond c\in \Sg^{\B}(X_1\cap X_2)$ and 
$(\Diamond a\land \square a_{1}\ldots \square a_{n}\to \Diamond c)$ and $(\Diamond c\to (\Diamond b_{1}\lor\ldots \Diamond b_{m}))$, so the
pair $(\Gamma ,\Delta )$ is separable. This is a contradiction.

\item  (ii) Suppose that $((a_{1}\land \ldots_{\B_n} a_{n})\to c)$ and $(c\to_{\B_n} b\lor b_{1}\lor\ldots b_{m}))$ for some $c\in \Sg^{\B}(X_1\cap X_2)$, 
where $\square a_{1},\ldots , a_{n}\in \Gamma $ and $\Diamond b_{1},\ldots, \Diamond b_{m}\in \Delta.$ 
Then $\square c\in \Sg^{\B}(X_1\cap X_2)$ and $((\square a_{1}\land\ldots  \square
a_{n})\to_{\B_n} \square c)$ and $(\square c\to_{\B_n} (\square b\lor \Diamond b_{1}\lor \ldots \Diamond b_{m}))$, so the pair $(\Gamma ,\Delta )$
is separable. Again this is a contradiction.
\endproof%
\end{description}Iterating the above claim and lemma \ref{neat} countably many types, we are done.
\end{demo}

We assume that $\B_n\subseteq \B_{\infty}$ for all $n\in \omega$.
Let $$W=\{(T,F)\in \B_{\infty}\times \B_{\infty}: (\exists n\in \omega) ((T,F) \text { is } (X_1, X_2, \B_n) \text { saturated })\}.$$
We can assume that $n$ is unique to $(T,F)$. Indeed for  $w=(T,F)\in W,$ set $D_w=\alpha+n\omega$ 
where $n$ is the least number such that $(T,F)$ is $(X_1, X_2, \B_n)$ saturated.

Now we define the accessability relation $R$ as follows:

If $(T,F)$ is $(X_1, X_2, \B_n)$ saturated and $(T'F')$ is $(X_1, X_2, \B_m)$ saturated (here $n$ and $m$ depends uniquely on 
$(T,F)$ and $(T', F')$ respectively,) then set  
$$(T,F)R(T'F')$$ iff $$n\leq m$$
and $$(\forall a\in \Sg^{\B_n}X_1)(\square a\in T\Rightarrow
a\in T^{\prime }),$$ 
and $$(\forall b\in \Sg^{\B_n}X_2)(\Diamond b\in F\Rightarrow
b\in F^{\prime }).$$
Let $$\mathfrak{K}=(W, R, D_{w})_{w\in W}.$$
We first define two maps on $\A_1=\Sg^{\A}X_1$ and $\A_2=\Sg^{\A}X_2$ respectively, then those will be pasted using the freeness of $\A$
to give the required single homomorphism. 
For a sets $I, M$ and $p\in {}^IM$, $^IM^{(p)}$ denotes the set 
$$\{x\in {}^IM:|\{i\in I: x_i\neq p_i\}|<\omega\}.$$
Now define $\psi_1: \Sg^{\A}X_1\to \mathfrak{F}_{\mathfrak{K}}$ by
$\psi_1(p)=(f_w)$ such that if $w=(T,F)\in \mathfrak{K}$ is $(X_1, X_2, \B_n)$ saturated,
and $D_w=\alpha+n\omega$, then for $x\in {}^{\alpha}D_w^{(Id)}$ 
$$f_w(x)=1\Longleftrightarrow {\sf s}_{x\cup (Id_{D_w\sim \alpha})}^{\B_n}p\in T.$$
Basically the extension of $x$ by the identity map, to avoid cumbersome possibly confusing notation, 
we shall denote the map $x\cup Id_{D_w\sim \alpha}$ simply by $x$, 
which is clearly a finite transformation.
Hence $$f_k(x)=1\Longleftrightarrow {\sf s}_{{x}}^{\B_n}p\in T.$$
Now define $\psi_2: \Sg^{\A}X_1\to \mathfrak{F}_{\mathfrak{K}}$ by
$\psi_1(p)=(f_w)$ such that if $w=(T,F)\in \mathfrak{K}$ is $(X_1, X_2, \B_n)$ saturated,
and $D_w=\alpha+n\omega$, then for $x\in {}^{\alpha}D_w^{(Id)}$ 
$$f_w(x)=1\Longleftrightarrow \neg {\sf s}_{x\cup (Id_{D_w\sim \alpha})}^{\B_n}p\in F.$$
Abusing notation slightly, we denote $\psi_1$ by $\psi$, and proceed to show that $\psi$ is a homomorphism. Then we show that $\psi_2$ is also a homorphism,
and that $\psi_1$ and $\psi_2$ agree on $\Sg^{\A}(X_1\cap X_2)$:

\begin{enumarab}

\item Let $\psi(p)=(g_w)$ and $\psi(q)=(f_w)$ and $\psi(p\lor q)=(h_w)$. We prove that for $x\in {}^{\alpha}D_w^{(Id)}$:
$$h_w(x)=1\Longleftrightarrow g_w(x)=1 \text { or }f_w(x)=1.$$
Assume that $h_w(x)=1$, then ${\sf s}_x(p\lor q)\in T$ or $\neg {\sf s}_x(p\lor q)\in F$. We have
\begin{equation*}
\begin{split}
&{\sf s}_x(p\lor q)\in T\\
&\Longleftrightarrow {\sf s}_xp\lor {\sf s}_xq\in T\\
&\Longleftrightarrow {\sf s}_xp\in T\text { or }{\sf s}_xq\in T\\
&\Longleftrightarrow g_w(x)=1\text { or } f_w(x)=1\\
\end{split}
\end{equation*}
(because $T$ is a theory and $(T,F)$ is saturated.)
\item  Let $\psi(p)=(f_w)$ and $\psi(\neg p)=(g_w)$. We show that $f_w(x)=0$ iff $g_w(x)=1$. We prove $\Rightarrow$.
Assume that $g_w(x)=1$.Then ${\sf s}_x(\neg p)\in T$ or $\neg({\sf s}_x\neg p)\in F$. We have
\begin{equation*}
\begin{split}
&{\sf s}_x(\neg p)\in T\\
&\Longleftrightarrow \neg {\sf s}_xp\in T\\
&\Longleftrightarrow {\sf s}_xp\notin T\text { and }\neg {\sf s}_xp\notin F\\
&\Longleftrightarrow f_k(x)=0\\
\end{split}
\end{equation*}

\item $\psi({\sf c}_kp)=(f_k)$. Assume that $f_k(x)=1$. Then 
${\sf s}_x{\sf c}_kp\in T$ or ${\sf -s}_x{\sf c}_kp\in F$. 
Assume the former. 
Let $$l\in \{\mu\in \beta: x^{-1}\{\mu\}=\{\mu\}\}\sim \Delta p.$$
Such an $l$ exists because $\B$ is dimension complemented.
Let $$\tau=x\circ [k,l].$$
Then $${\sf c}_l{\sf s}_{\tau}p={\sf s}_{\tau}{\sf c}_kp={\sf s}_x{\sf c}_kp,$$
and by the choice of $F,$ we have
$${\sf c}_l{\sf s}_{\tau}p\in F\Longleftrightarrow {\sf s}^l_u{\sf s}_{\tau}p\in F.$$
We use the following helpful notation. 
For a function $f$, the function $g=f(a\to u)$ is defined by $g(x)=f(x)$ for $x\neq a$ and $g(a)=u$. Now we have
\begin{equation*}
\begin{split}
\psi({\sf c}_kp)(x)=1
&\Longleftrightarrow {\sf s}_x{\sf c}_kp\in T\\
&\Longleftrightarrow {\sf c}_l{\sf s}_{\tau}p\in T\\
&\Longleftrightarrow {\sf s}_u^l{\sf s}_{\tau}p \in T\\
&\Longleftrightarrow {\sf s}_{x(k\to u)}p\in T\\
&\Longleftrightarrow x\in {\sf c}_k\psi(p)\\
\end{split}
\end{equation*}
\item Let $\psi(p)=(f_w)$ and $\psi(\square p)=(h_w)$ and $\square(\psi(p))=(g_w)$
$h_w(x)=1$. Let $w=(T,F).$ If  ${\sf s}_x\square p\in T$, then 
${\sf s}_xp\in T'$ for all $(T',F')$ such that $(T,F)R(T', F')$, hence $f_{w'}(x)=1$ for all $w'Rw$, then $g_w(x)=1$. 
If $$\neg {\sf s}_x\square p\in F.$$
then $\neg s_x \neg \diamond \neg p\in F$, then $s_x\diamond \neg p\in F$, then $s_x\neg p\in F'$ for all $F'$.
\end{enumarab}
\begin{theorem}
Let $\A$ be as in the above theorem, and $R$ be the accessibility relation constructed in the proof. Then, the following hold
\begin{table}[ht]\centering
\begin{tabular}{|c|c|c|c|c|}
\hline Formula & \multicolumn{2}{c|}{Condition on $R$}  \\ 
\hline $\Diamond T$ &  $\forall x\exists y Rxy$ & \\ 
\hline $\square p\to \square\square p$ & transitive: $xRy\land yRz\Rightarrow yRz$& \\ 
\hline $\square p\to p $& reflexive: $xRx$  &\\  
\hline
\end{tabular} 
\end{table}

\begin{enumarab}
\item If $\A\models \Diamond 1$ then $R$ satisfies $(\forall x)(\exists y xRy)$
\item if $\A\models \neg \square p\lor p=1$ then $R$ is reflexive
\item If $\A\models \neg \square p\lor \square \square p=1$, then $R$ is transitive.
\end{enumarab}
\end{theorem}

\begin{theorem} Let $\kappa$ be any ordinal $>1$. 
Let $\bold M=\{\A\in K_{\kappa+\omega}: \A=\Sg^{\A}\Nr_{\kappa}\A\}$.
Then $\bold M$ has $SUPAP$. 
\end{theorem}
\begin{demo}{Proof} We now first prove that $K$ has the amalgamation property. Let $\kappa$ be an arbitrary ordinal $>0.$
Let $\A,\B$ and $\C$ be in $K$ and $f:\C\to \A$ and $g:\C\to \B$ be monomorphisms. We want to find an amalgam.
Let $\langle a_i: i\in I\rangle$ be an enumeration of $\A$ and $\langle b_i:i\in J\rangle$ 
be an enumeration of $\B$ such that
$\langle c_i: i\in I\cap J\rangle$ is an enumeration of $\C$ with $f(c_i)=a_i$ and $g(c_i)=b_i$ for all $i\in I\cap J$. 
Then $\Delta a_i=\Delta b_i$ for all $i\in I\cap J$. Let $k=I\cup J$. 
Let $\xi$ be a bijection from $k$ onto a cardinal $\mu$.
Let $\rho\in {}^{\mu}\wp(\kappa+\omega)$ be defined by
$\rho \xi i=\Delta a_i$ for $i\in I$ and $\rho \xi j=\Delta b_j$ for $j\in J$. Then $\rho$ is well defined.
Let $\beta=\kappa+\omega$ and let $\Fr=\Fr_{\mu}^{\rho}K_{\beta}$.
Let $\Fr^{I}$ be the subalgebra of $\Fr$ generated by $\{\xi i/Cr_{\mu}^{\rho}K_{\beta} :i \in I\}$ and let 
$\Fr^J$ be the subalgebra generated by $\{\xi i/Cr_{\mu}^{\rho}K_{\beta} :i \in J\}.$
To avoid cumbersome notation we write $\xi i$ instead of $\xi i/Cr_{\mu}^{\rho}K_{\beta}$ and similarly for $\xi j$.
No confusion is likely to ensue. Then there exists 
a homomorphism from $\Fr^{I}$ onto $\A$
such that
$\xi i\mapsto a_i$ $(i\in I)$
and similarly 
a homomorphism from $\Fr^{J}$ into $\B$ such that
$\xi j\mapsto b_j$ $(j\in I).$
Therefore there exist ideals $M$ and $N$ ideals of $\Fr^I$ and $\Fr^J$ respectively, and 
there exist isomorphisms 
$$m:\Fr^{I}/M\to \A\text { and }n:\Fr^J/N\to \B$$ such that
$$m(\xi i/M)=a_i\text { and }(\xi i/N)=b_i.$$
Let $\Fr^{(I\cap J)}$ denote the subalgebra of $\Fr$ generated by $\{\xi i:i\in I\cap J\}$.
Then there is a homomorphism from $\Fr^{(I\cap J)}$ into $\Fr^{I}/M$ such that
$\xi i\mapsto \xi i/M$. 
Therefore
$$\theta:\Fr^{(I\cap J)}/M\cap \Fr^{(I\cap J)}\to (\Fr^I/M)^{(I\cap J)}$$ 
defined by 
$$\theta \xi i/M\cap \Fr^{(I\cap J)}=\xi i/M$$
is an isomorphism. Here 
$(\Fr^I/M)^{(I\cap J)}$ is the subalgebra of $\Fr^I/M$ generated by $\{\xi i/M:i\in I\cap J\}.$
But $$\psi: (\Fr^{I}/M)^{(I\cap J)}\to \C$$ defined by
$$\psi(\xi i/M)=c_i$$ 
is a well defined isomorphism.
Thus the map $\psi\circ \theta$ defined by 
$$\xi i/M\cap \Fr^{(I\cap J)}\mapsto c_i$$ is an isomorphism 
from $\Fr^{(I\cap J)}/M\cap \Fr^{(I\cap J)}$ onto $\C$.
Similarly the map $$\xi i/N\cap \Fr^{(I\cap J)}\mapsto c_i$$ 
is an isomorphism from  
$\Fr^{(I\cap J)}/N\cap \Fr^{(I\cap J)}$ onto $\C.$
It follows that $$M\cap \Fr^{(I\cap J)}=N\cap \Fr^{(I\cap J)}.$$
Now let $x\in \Ig(M\cup N)\cap \Fr^{I}$.
Then there exist $b\in M$ and $c\in N$ such that $x\leq b+c$. Thus $x-b\leq c$.
But $x-b\in \Fr^{(I)}$ and $c\in \Fr^{J}$, it follows that there exists an interpolant 
$d\in \Fr^{(I\cap J)}$ such that $x-b\leq d\leq c$. We have $d\in N$
therefore $d\in M$, and since $x\leq d+b$, therefore $x\in M$.
It follows that
$\Ig(M\cup N)\cap \Fr^{I}=M$
and similarly
$\Ig(M\cup N)\cap \Fr^{J}=N$.
In particular $P=\Ig(M\cup N)$ is a proper ideal.
Let $\D=\Fr/P$.
Let $k:\Fr^{I}/M\to \Fr/P$ be defined by $k(a/M)=a/P$
and $h:\Fr^{J}/M\to \Fr/P$ by $h(a/N)=a/P$. Then 
$k\circ m$ and $h\circ n$ are one to one and 
$k\circ m \circ f=h\circ n\circ g$. 

We now prove that $\Fr/P$ is actually a superamalgam.
i.e we prove that $K$ has the superamalgamation property.
Assume that $k\circ m(a)\leq h\circ n(b)$.
There exists $x\in \Fr^{I}$ such that $x/P=k(m(a))$ and $m(a)=x/M$.
Also there exists $z\in \Fr^{J}$ such that $z/P=h(n(b))$ and $n(b)=z/N$.
Now $x/P\leq z/P$ hence $x-z\in P$.
Therefore  there is an $r\in M$ and an $s\in N$ such that $x-r\leq z+s$.
Now $x-r\in \Fr^I$ and $z+s\in \Fr^J,$ it follows that there is an interpolant  $u\in \Fr^{(I\cap J)}$ such that 
$x-r\leq u\leq z+s$.
Let $t\in C$ such that $m\circ f(t)=u/M$ and $n\circ g(t)=u/N.$ 
We have  $x/P\leq u/P\leq z/P$.
Now $m(f(t))=u/M\geq x/M=m(a).$
Thus $f(t)\geq a$. Similarly 
$n(g(t))=u/N\leq z/N=n(b)$, 
hence $g(t)\leq b$. 
By total symmetry, we are done.
\end{demo}

\end{document}